\documentclass[12pt, normalheadings, DIV14]{scrartcl}  
  \usepackage{scrpage} 
  \usepackage{palatino, euler}
  \usepackage[notoday, nofancy]{rcsinfo}  

  \deftripstyle{mybodystyle}
   { Alexander and Thurston norms of fibered
      3-manifolds}{}{\rm \arabic{page}}
    {Version: \rcsInfoRevision}{}{\today}
  \newcommand{\firstpairscale}{1.0}
  \newcommand{\secondpairscale}{1.0}
  \newcommand{\linkscale}{0.9}
  \newcommand{\thirdpairscale}{1.0}

\usepackage{amsmath, amssymb, amsthm,  latexsym}
\usepackage{epsfig}

\usepackage[bf]{caption}

\newcommand{\euler}[1]{\chi(#1)}
\newcommand{\eulerminus}[1]{\chi_{-}(#1)}
\DeclareMathOperator{\ab}{ab}
\DeclareMathOperator{\Newt}{Newt}
\DeclareMathOperator{\Ann}{Ann}
\DeclareMathOperator{\Burau}{Burau}

\newcommand{\R}{{\mathbb R}}
\newcommand{\Z}{{\mathbb Z}}
\newcommand{\bdry}{\partial}
\newcommand{\iso}{\cong}
\newcommand{\cross}{\times}
\newcommand{\tensor}{\otimes}
\newcommand{\directsum}{\oplus}
\newcommand{\maps}{\colon\thinspace}
\newcommand{\Hom}{ {\mathrm{Hom}}}
\newcommand{\GL}{ {\mathrm{GL}} }
\DeclareMathOperator{\rank}{rank}
\DeclareMathOperator{\Aut}{Aut}
\newcommand{\norm}[1]{{\left\|  #1  \right\|}}
\newcommand{\class}[1]{{\left[ \, #1  \, \right]}}
\newcommand{\abs}[1]{{\left| #1 \right|}}

\newcommand{\setdef}[2]{{  \left\{  {#1}  \ \left| \   {#2} \right. \right\} }}
\newcommand{\spandef}[2]{{  \left\langle  {#1}  \ \left| \   {#2} \right. \right\rangle }}
\newcommand{\mtext}[1]{\quad\mbox{#1}\quad}

\newtheoremstyle{plain}{}{}{\slshape}{}{\bfseries}{.}{0.5em}{}
\theoremstyle{plain} 

\swapnumbers
\newtheorem{theorem}{Theorem}[section]

\newtheorem{lemma}[theorem]{Lemma}

\newtheorem{proposition}[theorem]{Proposition}

\makeatletter
  \let\c@theorem=\c@subsection
  \let\c@figure=\c@subsection
  \let\p@figure=\p@subsection
  
  \let\cl@figure=\cl@subsection
\makeatother

\newcommand{\McMullensQ}{Question~A}
\newcommand{\MyQ}{Question~B}
\newcommand{\MyModQ}{Question~$\text{B}'$}
\newtheorem*{McMullensQuestion}{\McMullensQ}
\newtheorem*{MyQuestion}{\MyQ}
\newtheorem*{MyModQuestion}{\MyModQ}
\newtheorem*{answer}{Answer}

\begin{document}

   \rcsInfo $Id: norms_of_fibered.tex,v 1.3 2000/07/13 18:58:15 nathand Exp $

   \thispagestyle{plain}
   \

   \vspace{1in}
   \noindent
   {\sf \large \textbf{Alexander and Thurston norms of fibered 3-manifolds}}

   \bigskip

   \noindent
   {Nathan M.~Dunfield\footnote{This work was partially supported by a
    Sloan Dissertation Fellowship.}}

   \medskip

   \noindent
   {\small Department of Mathematics, Harvard University, One Oxford St.,
    Cambridge, MA 02138\\
    email:  nathand@math.harvard.edu}

   \bigskip
 
   \noindent
   {\small Version: \rcsInfoRevision, Compile: \today, Last commit:
    \rcsInfoLongDate}

\section{Introduction}

\subsection{Statement of results}\label{StatementOfResults}

For a 3-manifold $M$, McMullen derived from the Alexander polynomial
of $M$ a norm on $H^1(M, \R)$ called the Alexander norm.  He showed
that the Thurston norm on $H^1(M, \R)$, which measures the complexity
of a dual surface, is an upper bound for the Alexander norm.  He asked
(\McMullensQ\ below) if these two norms were equal on all of
$H^1(M,\R)$ when $M$ fibers over the circle.  Here, I will give
examples which show that the answer to \McMullensQ\ is emphatically
no.  As explained below, \McMullensQ\ is related to the faithfulness
of the Gassner representations of the braid groups.  The key tool used
to understand \McMullensQ\ is the Bieri-Neumann-Strebel invariant from
combinatorial group theory.  Theorem~\ref{B_AandSigma_A} below, which is of
independent interest, connects the Alexander polynomial with a certain
Bieri-Neumann-Strebel invariant.

I will begin by reviewing the definitions of the Alexander and
Thurston norms, and Theorem~\ref{NormThm} which relates them.
Then I'll discuss \McMullensQ\ and the connection to the braid groups.
After that, I'll state \MyQ, a much weaker version of \McMullensQ, to
which the answer is also no.  A brief description of the examples
which answer these two questions concludes
Section~\ref{StatementOfResults}.  In Section~\ref{MoralReason}, I'll
connect these questions with the Bieri-Neumann-Strebel invariants, and
explain why, morally speaking, the answer to both questions must be
no.  Section~\ref{DescriptionOfRest} outlines the rest of the paper.

The Alexander norm is defined in \cite{McMullenNorm} as follows.  Let
$M$ be a 3-manifold (all 3-manifolds in this paper will be assumed to
be connected).  Let $G$ be the fundamental group of $M$.  Let $\ab(G)$
denote the maximal free abelian quotient of $G$, which is isomorphic
to $\Z^{b_1(M)}$ where $b_1(M)$ is the first Betti number of $M$.  The
Alexander polynomial $\Delta_M$ of $M$ is an element of the group ring
$\Z[\ab(G)]$.  It is an invariant of the homology of the cover of $M$
with covering group $\ab(G)$ (for details see
Section~\ref{DefOfAlexPoly}).  The Alexander norm on $H^1(M, \R)$ is
the norm dual to the Newton polytope of $\Delta_M$.  That is, if
$\Delta_M = \sum_{i=1}^n a_i g_i$ with $a_i \in \Z \setminus \left\{ 0
\right\}$ and $g_i \in \ab(G)$ then the norm of a class $\phi \in
H^1(M,\R)$ is defined to be
\[
\norm{\phi}_A = \sup_{i,j} \phi(g_i - g_j).
\]
The unit ball $B_A$ of this norm is, up to scaling, the polytope dual to
the Newton polytope of $\Delta_M$.  

The Thurston norm is defined as follows.  For a compact connected
surface $S$, let $\eulerminus{S} = \abs{\euler{S}}$ if $\euler{S} \leq
0$ and $0$ otherwise.  For a surface with multiple connected
components $S_1, S_2, \dots, S_n$, let $\eulerminus{S}$ be sum of the
$\eulerminus{S_i}$.  Then the Thurston norm of an integer class $\phi
\in H^1(M, \Z) \iso H_2(M, \bdry M; \Z)$ is
\begin{equation*}
\begin{split}
\norm{\phi}_T = \inf \left\{ \eulerminus{S} \ \right| &\  \mbox{$S$ is a properly embedded oriented surface}\\
  &\left. \mbox{that is dual to $\phi$} \right\}.
\end{split}
\end{equation*}
As described in \cite{ThurstonNorm}, this norm extends continuously to
all of $H^1(M, \R)$.  The unit ball $B_T$ in this norm is a
finite-sided convex polytope.
  
It should be noted that both of these ``norms'' are really
semi-norms---they can be zero on non-zero vectors of $H^1(M,\R)$.  

McMullen proved the following theorem which connects the two norms;
here $b_i(M) = \rank H_i(M, \R)$ denotes the $i^{\text{th}}$ Betti number of
$M$.
\begin{theorem}[\cite{McMullenNorm}]\label{NormThm}
Let $M$ be a compact, orientable 3-man\-i\-fold whose boundary,
if any, is a union of tori.   Then for all
$\phi$ in $H^1(M, \R)$, the Alexander and Thurston norms satisfy
\[  
\norm{\phi}_A \leq \norm{\phi}_T \mtext{if $b_1(M) \geq 2,$}
\]
or
\[
\norm{\phi}_A \leq \norm{\phi}_T + 1 + b_3(M) \mtext{if $b_1(M) = 1$
  and $\phi$ generates $H^1(M, \Z)$.}
\]
Moreover, equality holds when $\phi \maps \pi_1(M) \to \Z$ and $\phi$
can be represented by a fibration $M \to S^1$, where the fibers have
non-positive Euler characteristic.
\end{theorem}
This theorem generalizes the fact that the degree of the Alexander
polynomial of a knot is bounded by twice the genus of any Seifert
surface.  In many simple cases, e.g.~almost all the exteriors of the
links with $9$ or fewer crossings, the Alexander and Thurston norms
agree on all of $H^1(M, \R)$ (see \cite{McMullenNorm}).  In such
cases, this theorem explains D.~Fried's observation from the 80's that
frequently the shape of the Newton polytope of the Alexander
polynomial is dual to that of the Thurston norm ball.

Before stating \McMullensQ, I need to discuss the relationship between
the Thurston norm and cohomology classes $\phi \maps \pi_1(M) \to \Z$
which can be represented by fibrations $M \to S^1$.  There are
top-dimensional faces, called the \emph{fibered faces}, of $B_T$ such
that a class $\phi \in H^1(M, \Z)$ can be represented by a fibration
over the circle if and only if $\phi$ lies in the cone over the
interior of one of the fibered faces \cite[$\S 3$]{ThurstonNorm}.  In
this context, the last sentence of Theorem~\ref{NormThm} is equivalent
to ``Moreover, the two norms agree on classes that lie in the cone
over the fibered faces of $B_T$''.  The point of this paper is to
answer:
\begin{McMullensQuestion}[McMullen \cite{McMullenNorm}]
Let $M$ be a compact, orientable 3-man\-i\-fold whose boundary,
if any, is a union of tori.  Suppose that $M$ fibers over the circle and
that $b_1(M) \geq 2$.  Do the Alexander and Thurston norms agree on all
of $H^1(M, \R)?$
\end{McMullensQuestion}

My motivation for studying this question is McMullen's result that a
yes answer would imply that the Gassner representations of the pure
braid groups are all faithful \cite[$\S 8$]{McMullenNorm}.  This would
answer in the affirmative the important question: Are the braid groups
linear, that is, do they have faithful, finite-dimensional, linear
representations?  Sadly, I will show that the answer to \McMullensQ\ 
is no in a strong sense.  (Note: Since I wrote this paper, Bigelow and
Krammer have independently shown that braid groups are linear
\cite{Bigelow00, Krammer99, Krammer00}.  Their proofs use a different
representation, and it remains unknown whether the Gassner
representation is faithful).

\begin{figure}[htbp]
\centering
\parbox[t]{0.3\textwidth}{\centering
    \psfig{figure=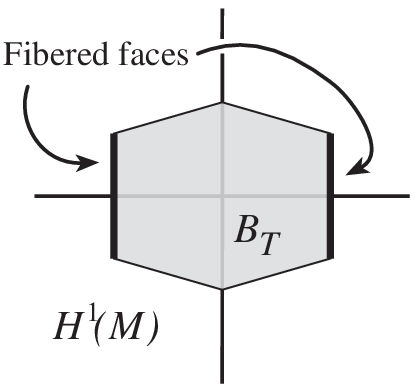, scale=\firstpairscale}
    \caption{ The Thurston norm ball.}\label{B_T}
    }\hspace{0.05\textwidth}%
\parbox[t]{0.6\textwidth}{\centering 
     \psfig{figure=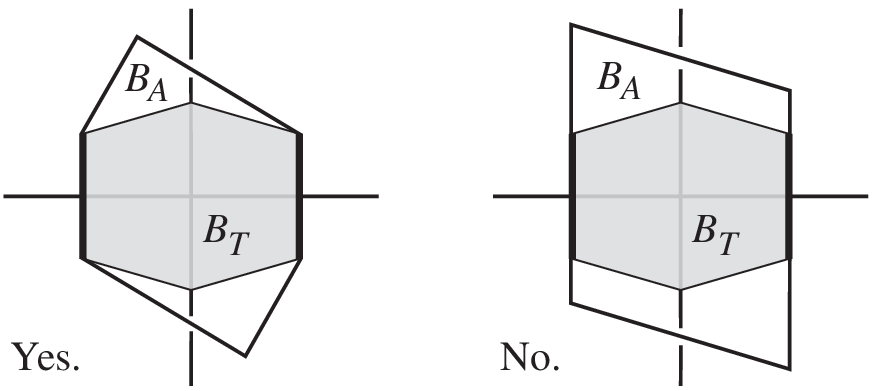, scale=\firstpairscale}
     \caption{Possible answers to \MyQ.}\label{ComparingF_AandF_T}
     }
\end{figure}

To explain why the answer to \McMullensQ\ is no, let me formulate a
weaker version of \McMullensQ\ which will help make clear some of the
issues involved.  Henceforth, I will assume that $b_1(M) \geq 2$.  A
typical example of $B_T$ is given in Figure~\ref{B_T}.

  There is a pair of fibered faces and the rest of
the faces are not fibered.  Theorem~\ref{NormThm} tells us that
$\norm{ \ \cdot \ }_A \leq \norm{\ \cdot \ }_T$ hence that $B_A
\supseteq B_T$.  Since the two norms agree on a fibered face $F_T$ of
$B_T$, there is a face $F_A$ of $B_A$ which contains $F_T$.  Now, it
seems a bit much to expect that if $M$ fibers over the circle then the
two norms agree on classes that are far from any fibered face.  So
it's reasonable to consider:
\begin{MyQuestion}  
Let $M$ be a compact, orientable 3-manifold whose boundary,
if any, is a union of tori.  Suppose that $M$ fibers over the circle
and that $b_1(M) \geq 2$.  Let $F_T$ be a fibered face of $B_T$ and
$F_A$ the face of $B_A$ which contains it.  Are $F_T$ and $F_A$ always
equal?
\end{MyQuestion}
Figure~\ref{ComparingF_AandF_T} shows the two possibilities.  Note
that a yes answer to \McMullensQ\ implies a yes answer to \MyQ.  I
will give examples which show that
\begin{answer} The answer to \MyQ, and therefore
\McMullensQ, is no.
\end{answer}
I will give two kinds of examples.  In Section~\ref{BraidGroups}, I
will constuct examples using the fact that the Burau representation of
the braid group on 5 strands is not faithful.
Section~\ref{BraidGroups} is independent of the rest of the paper.
Section~\ref{SpecificExample} contains an example which is the
exterior of a specific $17$ crossing link in $S^3$.

McMullen's formulation of \McMullensQ\ restricted attention to those
manifolds which are the exteriors of links in $S^3$.  All my examples
are such manifolds, but I felt the more general statement was
appropriate here.

\subsection{Connection to the BNS invariants}\label{MoralReason}

In this section I will describe the connection between \MyQ\ and the
Bieri-Neumann-Strebel (BNS) invariants.  In light of this connection,
I will explain why the answer to \MyQ\ must be, morally speaking, no.
The BNS invariants will also be used in constructing and verifying the
example in Section~\ref{SpecificExample}.

I'll begin with the definition of the BNS invariants (for details see
\cite{BieriNeumannStrebel87}, and from a different point of view,
\cite{Brown87}).  Let $G$ be a finitely-generated group.  Set
\[
S(G) = \left( H^1(G, \R)  \setminus \{ 0 \} \right) \big/ \ \R^+, 
\]
where $\R^+$ acts by scalar multiplication and $S(G)$ is given the
quotient topology.  A point $\class{\chi}$ in $S(G)$ will be thought of as
an equivalence class of homomorphisms $\chi \maps G \to \R$.  For $\class{\chi} \in
S(G)$ define  $G_\chi = \chi^{-1}\left([0, \infty)\right) = \setdef{g
  \in G}{\chi(g) \geq 0}$, which is a sub-monoid of $G$.  

Let $H$ be a group acted on by $G$ where $G'$ acts by inner
automorphisms (e.g. $H = G'$ where $G$ acts by conjugation).  Then the
BNS invariant of $G$ and $H$ is:
\[
\begin{split}
\Sigma_H = \left\{ \class{\chi} \in S(G) \ \right| &\ \mbox{$H$
 is finitely generated over some}\\
   &\left. \mbox{finitely generated sub-monoid of $G_\chi$} \right\}.
\end{split}
\]
It turns out that $\Sigma_H$ is always an open subset of the sphere
$S(G)$.  

Let $M$ be a 3-manifold, and $G = \pi_1(M)$.   Set $\Sigma= \Sigma_{G'}$.  Bieri,
Neumann, and Strebel proved the following with the help of Stallings'
fibration theorem:
\begin{theorem}[{\cite[Thm.~E]{BieriNeumannStrebel87}}]
Let $M$ be a compact, orientable, irreducible 3-manifold.  Then
$\Sigma$ is exactly the projection to $S(G)$ of the interiors of the
fibered faces of the Thurston norm ball $B_T$.
\end{theorem}

For convenience, in the rest of this section I will assume that
$H_1(M, \Z)$ is free.  This is not essential, and the theory will be
developed without this assumption in
Sections~\ref{AlexanderPoly}-\ref{Alex_and_BNS}.  The commutator
subgroup $G'$ is the fundamental group of the universal abelian cover
of $M$.  So $A = G'/G''$ is the first homology of that cover.  Thought
of as a module over $\Z[\ab(G)]$, $A$ is the Alexander invariant of
$M$, from which the Alexander polynomial is derived.  Thus it is not
too surprising that the BNS invariant $\Sigma_A$ is connected to the
Alexander polynomial:
\begin{theorem} \label{B_AandSigma_A}
Let $M$ be a compact, orientable 3-manifold.  There are
top-dimensional faces $F_1, F_2, \dots, F_n$ of the Alexander ball
$B_A$ such that the projection of the interiors of the $F_i$ into
$S(G)$ is exactly $\Sigma_A$.  Moreover, the $F_i$ are completely
determined by the Alexander polynomial of $M$.
\end{theorem}
Theorem~\ref{fullB_AandSigma_A} below is an expanded version of
Theorem~\ref{B_AandSigma_A} which explains how the $F_i$ are
determined.  Now since $A$ is a quotient of $G'$, it follows
immediately from the definitions that $\Sigma_A \supset \Sigma$.
Combining this with Theorem~\ref{B_AandSigma_A}, it follows that \MyQ\ 
is equivalent to:
\begin{MyModQuestion}
Let $M$ be a compact,  orientable 3-manifold whose boundary,
if any, is a union of tori.  Suppose that $M$ fibers over the circle
and that $b_1(M) \geq 2$.  Let $C$ be a connected component of
$\Sigma$.  If $D$ is the connected component of $\Sigma_A$ which
contains $C$, is $D$ always equal to $C$?
\end{MyModQuestion}
Put this way it begins to become clear that the answer to \MyQ\
should be no.  For many groups $G$, $\Sigma_{G'}$ is strictly contained
in $\Sigma_{G'/G''}$.  It remains only to produce examples of
3-manifolds whose fundamental groups have this property.  

\subsection{Outline of rest of paper}\label{DescriptionOfRest}
Section~\ref{BraidGroups} describes how to construct examples using
the Burau representation.  Section~\ref{AlexanderPoly} defines the
Alexander polynomial and proves a fact about the Alexander invariant
that's needed to prove Theorem~\ref{B_AandSigma_A}.
Section~\ref{BNS_definitions} discusses the BNS invariants and records
the properties that will be needed later.  Section~\ref{Alex_and_BNS}
proves the full version of Theorem~\ref{B_AandSigma_A}.  Finally,
Section~\ref{SpecificExample} gives an example of a specific link
exterior in $S^3$ for which the answer to \MyQ\ is no.

\subsection{Acknowledgments}  I wrote this
paper while a graduate student at the University of Chicago supported
by a Sloan Dissertation Fellowship.  I would like to thank Curt
McMullen for useful correspondence.  I got interested in the
connection between the Alexander polynomial and Thurston norm at a
problem session at KirbyFest (MSRI, June 1998), where Fried's
observation was related by Joe Christy.  I would like to thank the
organizers, MSRI, and the NSF for support to attend that conference.
I would also like to thank the referee for help clarifying the
exposition.

\section{Connection with braid groups}\label{BraidGroups}

Let $B_n$ denote the $n$-strand braid group.  McMullen showed that if
the answer to \McMullensQ\ is yes,  then the Gassner representation of
$B_n$ is faithful for all $n$ \cite{McMullenNorm}.
In this section, I'll give a very similar argument to show: 
\begin{proposition}\label{BurauAndQuestions} If the answer to \MyQ\ is yes, then the Burau
representation of $B_n$ is faithful for all $n$.
\end{proposition}
Since the Burau representation of $B_n$ is \emph{not} faithful for $n
\geq 5$ \cite{Bigelow99, LongPaton, Moody93}, the proposition implies
that the answer to \MyQ, and hence \McMullensQ, is no.  

Before proving the proposition, let me define the braid groups and the
Burau representation (see  \cite{Birman74} for more).  Let $D_n$ be the
disc with $n$ punctures.  Consider the group of homeomorphisms
$\Hom^+(D_n, \bdry D_n)$ of $D_n$ which are orientation preserving and
fix $\bdry D_n$ pointwise.  The braid group $B_n$ is
$\Hom^+(D_n, \bdry D_n)$ modulo isotopies which pointwise fix $\bdry D_n$.

To define the Burau representation, consider the homomorphism
\[
\phi\maps H_1(D_n) \to \Z = \left< t \right>
\]
which takes any clockwise oriented loop about a single puncture to
$t$.  Let $\tilde{D}_n$ be the cover of $D_n$ corresponding to $\phi$.
The homology of $\tilde{D}_n$ is a module over the group ring
$\Z\left[\left< t \right> \right]$ of the group of covering
transformations.  The module $H_1(\tilde{D}_n, \Z)$ is free of rank
$n-1$.  The Burau representation is a homomorphism $\Burau \maps B_n
\to \Aut \big( H_1 (\tilde{D}_n) \big)$.  By
$\Aut\big(H_1(\tilde{D}_n)\big)$, I mean automorphisms of
$H_1(\tilde{D}_n)$ as a $\Z\left[\left< t \right> \right]$-module.
Choosing a $\Z\left[\left< t \right> \right]$ basis of
$H_1(\tilde{D}_n)$ allows one to view the Burau representation as
having image in $\GL(n-1, \Z\left[\left< t \right> \right])$.  Given
$\beta$ in $B_n$, $\Burau(\beta)$ is constructed as follows.  Let $f
\maps D_n \to D_n$ be a representative of $\beta$.  Choose a lift
$\tilde{f} \maps \tilde{D}_n \to \tilde{D}_n$ of $f$.  Since the
action of $f$ on $H_1(D_n)$ commutes with $\phi$, the lift $\tilde{f}$
is equivariant.  Thus there is a unique lift of $f$ which leaves the
inverse image of $\bdry D_n$ pointwise fixed.  Let $\tilde{f}$ be that
lift and set $\Burau(\beta) = \tilde{f}_* \maps H_1(\tilde{D}_n) \to
H_1(\tilde{D}_n)$.

I'll need the following property of the Burau representation (see also
\cite{Morton99}).  Suppose $\beta$ is a braid whose action on the set
of punctures is an $n$-cycle.  Let $M_\beta$ be the 3-manifold which
is the mapping torus of $\beta$.  The manifold $M_\beta$ has two
boundary components, and $H_1(M_\beta) = \Z \directsum \Z$.  Take as a
basis of $H_1(M_\beta)$ the pair $(t', w)$ where $t'$ is a
counter-clockwise loop about a puncture in $D_n$ and $w$ is a point in
$\bdry D_n$ cross $S^1$.  It's not hard to see that the universal
abelian cover of $M_\beta$ is $\tilde{D}_n \cross \R$.  The covering
transformation corresponding to $t'$ is $\big(\tilde{d}, r\big)
\mapsto \big(t(\tilde{d}), r\big)$, and the covering transformation
corresponding to $w$ is $\big(\tilde{d}, r\big) \mapsto
\big(\tilde{f}(\tilde{d}), r + 1\big)$.  If we replace $t$ by $t'$ in
$\Burau(\beta)$, the matrix $( w I - \Burau(\beta))$ is a presentation
matrix for the homology of the universal abelian cover of $M_\beta$ as
a $\Z[H_1(M_\beta)]$-module.  Thus
\[
\Delta_{M_\beta} = \det\big( w I - \Burau(\beta)\big).
\]

I will now prove the proposition.

\begin{proof}[Proof of Proposition~\ref{BurauAndQuestions}]
Suppose the answer to \MyQ\ is yes and the Burau representation of
$B_n$ has kernel for some $n$.  As the Burau representation is known
to be faithful for $n = 2$, assume $n$ is at least $3$.  Then there is
a pseudo-Anosov element $\delta$ in the kernel \cite{Long86,
  Ivanov92}.  Replacing $\delta$ with a power of $\delta$ if
necessary, we can assume $\delta$ is a pure braid, that is, fixes each
puncture.  Let $\gamma$ be the braid $\sigma_1 \sigma_2 \dots
\sigma_{n-1}$ where $\sigma_i$ is the $i^{\text{th}}$ standard
generator of $B_n$ (see Figure~\ref{braid_gamma}).
\begin{figure}[htbp]
  \centering 
  \parbox[b]{0.3\textwidth}{
    \centering
    \psfig{figure=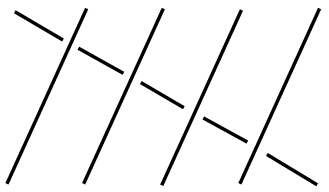, scale=\secondpairscale}
    }\hspace{0.05\textwidth}%
  \parbox[b]{0.6\textwidth}{
    \centering
    \psfig{figure=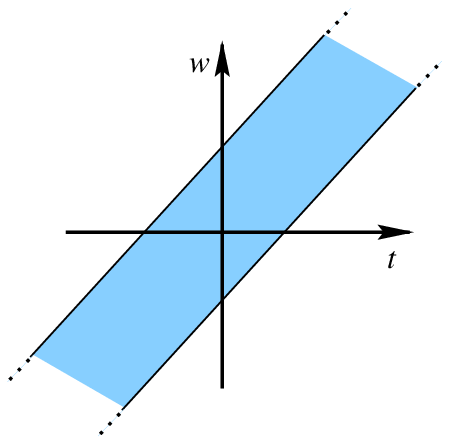, scale=\secondpairscale}
   }
  \parbox[t]{0.3\textwidth}{
    \centering
     \caption{The braid $\gamma$ when $n=5$.}\label{braid_gamma}
    }\hspace{0.05\textwidth}%
  \parbox[t]{0.6\textwidth}{
    \centering
    \caption{The Thurston norm ball of $M_\gamma$.}\label{norm_of_gamma}
    }
\end{figure}
Taking a power of $\delta$ if necessary, we can assume that $\beta =
\delta \gamma$ is pseudo-Anosov.  Now $\beta$ induces an $n$-cycle on
the punctures because $\delta$ was a pure braid and $\gamma$ induces
an $n$-cycle.  Since $\Burau(\beta) = \Burau(\gamma)$, the Alexander
polynomials of $M_\beta$ and $M_\gamma$ are the same.  The manifold
$M_\gamma$ is Seifert fibered, and it's easy to see that the Thurston
norm ball is as shown in Figure~\ref{norm_of_gamma}, where the two
infinite faces are fibered faces.  Thus by Theorem~\ref{NormThm}, the
Alexander norm ball of $M_\gamma$ has exactly the same shape as the
Thurston norm ball.  Since $M_\gamma$ and $M_\beta$ have the same
Alexander polynomials, the Alexander norm ball of $M_\beta$ is as
shown.  But $M_\beta$ is hyperbolic, and hence the Thurston norm is
non-degenerate.  So any face of the Thurston norm ball is bounded.
Thus a fibered face of the Thurston norm ball of $M_\beta$ is properly
contained in the corresponding face of the Alexander norm ball.  This
contradicts the assumption that the answer to \MyQ\ is yes.
\end{proof}

\section{The Alexander polynomial and its friends}\label{AlexanderPoly}

\subsection{Definitions}\label{DefOfAlexPoly}

I'll begin by reviewing the definition of the Alexander polynomial
and related invariants (for more see \cite{Hillman81, Rolfsen76,
  McMullenNorm}).  Let $X$ be a finite CW-complex with fundamental
group $G$.  Let $\tilde{X}$ be the universal free abelian cover of
$X$, that is, the cover induced by the homomorphism from $G$ to its
free abelianization $\ab(G)$.  Let $p$ be a point of $X$, and
$\tilde{p}$ its inverse image in $\tilde{X}$.  The
\emph{Alexander module} of $X$ is
\[
A_X = H_1(\tilde{X}, \tilde{p}; \Z)
\]
thought of as a module over the group ring $\Z[\ab(G)]$.  The reason
one uses the \emph{free} abelianization is so that the ring
$\Z[\ab(G)]$ has no zero divisors.

For a finitely generated module $M$ over $\Z[\ab(G)]$, the
$i^{\text{th}}$  elementary ideal $E_i(M) \subset \Z[\ab(G)]$ is
defined as follows.  Take any presentation
\[
0 \to \left( \Z[\ab(G)] \right)^r \xrightarrow{P}  \left( \Z[\ab(G)]
\right)^s \to M \to 0
\]
and set $E_i(M)$ to be the ideal generated by the $(s - i, s- i)$
minors of the matrix $P$.  The \emph{Alexander ideal} of $X$ is
$E_1(A_X)$.  The \emph{Alexander polynomial} of $X$, denoted
$\Delta_X$, is the greatest common divisor of the elements of the
Alexander ideal.  The polynomial $\Delta_X$ is defined up to
multiplication by a unit $g \in \ab(G)$ of $\Z[\ab(G)]$.
Equivalently, $\Delta_X$ is a generator of the smallest principle
ideal containing the Alexander ideal.

I should mention that the Alexander module, and hence Alexander
polynomial, depends only on the fundamental group of $X$;  it can be
thought of as an invariant of a finitely generated group.

I will need to consider $B_X = H_1(\tilde{X}; \Z)$, the
\emph{Alexander invariant} of $X$.  When $H_1(X; \Z)$ is free, $B_X =
G'/G''$.  As with $A_X$, the Alexander invariant $B_X$ is to be
thought of as a module over $\Z[\ab(G)]$.  The two modules are related
as follows.  Let $m \subset \Z[\ab(G)]$ be the augmentation ideal,
that is $m = \spandef{1 - g}{g \in \ab(G)}$.  The homology long exact
sequence for the pair $(\tilde{X}, \tilde{p})$ gives rise to the
short exact sequence
\[
0 \to B_X \to A_X \to m \to 0.
\]
The Alexander polynomial of $X$ could just have well been defined as
the gcd of $E_0(B_X)$ (for the equivalence of these two definitions see,
e.g. \cite{Traldi82}).

\subsection{Structure of the Alexander invariant of a 3-manifold}

The following fact about the structure of the Alexander ideal of a
3-manifold was crucial in McMullen's proof of
Theorem~\ref{NormThm}.
\begin{theorem}[{\cite[5.1]{McMullenNorm}}]\label{AlexIdealThm}
Let $M$ be a compact, orientable 3-man\-i\-fold whose boundary, if any, is
a union of tori.  Let $G = \pi_1(M)$.  Then $E_1(A_M) = m^p \cdot
(\Delta_M)$ where
\[
p = \begin{cases}
  0 & \text{if $b_1(M) \leq 1$, }\\
  1 + b_3(M) & \text{otherwise,}
\end{cases}
\]
and $m$ is the augmentation ideal of $\Z[\ab(G)]$.
\end{theorem}

The corresponding fact about $E_0(B_M)$ will be key to the proof of
Theorem~\ref{B_AandSigma_A}.  For a manifold with non-empty torus
boundary, Crowell and Strauss \cite{CrowellStrauss} showed that
$E_0(B_M) =(\Delta_M) \cdot m^q$ for an explicit value of $q$.  The
following proposition is weaker than \cite{CrowellStrauss}, but it
also applies to closed 3-manifolds.  It will suffice for my purposes and
follows easily from known results.
\begin{proposition}\label{Alex_prop}
Let $M$ be a compact, orientable 3-manifold whose boundary, if any, is
a union of tori.  Then
\[
\sqrt{E_0(B_M)} \cap m = \sqrt{ (\Delta_M ) } \cap m.
\]
\end{proposition}

\begin{proof}
By Theorem~1.1 of \cite{Traldi82} the short exact sequence
\[
0 \to B_M \to A_M \to m \to 0
\]
implies that there are integers $r, s \geq 0$ such that 
\[
E_1(A_M) \cdot m^r \subset E_0(B_M) \mtext{and} 
E_0(B_M)  \cdot m^s \subset E_1(A_M).
\]
Combining and multiplying by $m$ gives 
\[
E_1(A_M) \cdot m^{r + s + 1} \subset E_0(B_M) \cdot m^{s + 1} \subset
E_1(A_M) \cdot m.
\]
Taking radicals of the above and using that $\sqrt{I \cdot J} =
\sqrt{I} \cap \sqrt{J}$ gives 
\[
\sqrt{E_1(A_M)} \cap \sqrt{m} = \sqrt{E_0(B_M)} \cap \sqrt{m}.
\] 
Now $m$ is radical since it is the kernel of the ring homomorphism
$Z[\ab(G)] \to \Z$ which sends every $g \in \ab(G)$ to 1.  By
Theorem~\ref{AlexIdealThm} we have $E_1(A_M) = (\Delta_M) \cdot m^p$.
Combining, we get $\sqrt{E_0(B_M)} \cap m = \sqrt{ (\Delta_M ) } \cap
m$ as desired.
\end{proof}

\section{Bieri-Neumann-Strebel Invariants}
\label{BNS_definitions}

Recall the definition of the BNS invariant from
Section~\ref{MoralReason}.  Let $G$ be a finitely-generated group.
Let $S(G) = \left( H^1(G, \R) \setminus \{ 0 \} \right) \big/ \ \R^+
$.  For $\class{\chi} \in S(G)$ we have the sub-monoid $G_\chi = \setdef{g
  \in G}{\chi(g) \geq 0}$.  Let $H$ be a group acted on by $G$ where
$G'$ acts by inner automorphisms.  Then the BNS invariant of $G$ and $H$
is:
\[
\begin{split}
\Sigma_H = \left\{ \class{\chi} \in S(G) \ \right| &\ \mbox{$H$
 is finitely generated over some}\\
   &\left. \mbox{finitely generated sub-monoid of $G_\chi$} \right\}.
\end{split}
\]
We can also consider the larger invariant
\[
\Sigma'_H = \setdef{\class{\chi} \in S(G)}{\mbox{$H$ is finitely generated
    over $G_\chi$}}.
\]
When $H$ is abelian $\Sigma'_H = \Sigma_H$
\cite[Theorem~2.4]{BieriNeumannStrebel87}.  The special case of
$\Sigma'_H$ when both $G$ and $H$ are abelian was studied by Bieri and
Strebel \cite{BieriStrebel81} prior to the development of the full BNS
invariant.  The rest of this section will be devoted to that special
case.

Let $Q$ be a finitely generated free abelian group and $A$ a finitely
generated $\Z[Q]$-module.  Since $A$ has an action of $Q$, we can form
the BNS invariant $\Sigma_A = \Sigma'_A$.  To reduce clutter, I'll denote
$\Z[Q]$ by $\Z Q$.  A basic property shown in
\cite[\S 1.3]{BieriStrebel81} is that $\Sigma_A = \Sigma_{\Z Q/\Ann(A)}$
where $\Ann(A)$ is the annihilator ideal of $A$.  Thus $\Sigma$ can be seen
as an invariant of an ideal $I \subset \Z Q$.  The following basic
identities hold for any ideals $I, J$ in $\Z Q$
\cite[\S 1.3]{BieriStrebel81}:
\[
\Sigma_{\Z Q / I} = \Sigma_{\Z Q/\sqrt{I}} \mtext{and}
\Sigma_{\Z Q / (I \cdot J)} = \Sigma_{\Z Q/(I \cap J)} = 
\Sigma_{\Z Q/I} \cap \Sigma_{\Z Q/J} .
\]

For principle ideals $I$, the invariant $\Sigma_{\Z Q/I}$ can be
easily calculated, as I will now describe.  For $p \in \Z Q$, the
\emph{Newton polytope} $\Newt(p)$ is defined as follows.  Consider the
vector space $V = Q \tensor \R$ which contains $Q$ as a lattice.  The
Newton polytope of $p$ is the convex hull in $V$ of those $q \in Q$
which have non-zero coefficient in $p$.  The vertices of $\Newt(p)$ lie
in $Q$, and I'll define the coefficient of a vertex of $\Newt(p)$ to be
the non-zero coefficient of the corresponding term of $p$.  Given a
$q$ in $Q$, define the open hemisphere $H_q$ of $S(Q)$ to be
\[
\setdef{\class{\chi} \in  S(Q)}{\chi(q) >  0}.
\]
The following theorem allows us to calculate $\Sigma_{\Z Q/I}$ for a
principle ideal $I$.  

\begin{theorem}[{\cite[5.2]{BieriStrebel81}}]\label{computing_Sigma}
Let $Q$ be a finitely generated free abelian group and $p$ an element of
$\Z Q$.  The connected components of $\Sigma_{\Z Q/(p)}$
are in one-to-one correspondence with the vertices of $\Newt(p)$ whose
coefficients are $\pm 1$, where such a vertex $v$ corresponds to
\[
C_v = \bigcap \setdef{H_{vw^{-1}}}{\text{$w$ is a vertex of $\Newt(p)$
      distinct from $v$}}.
\]
\end{theorem}

\section{BNS invariants and Alexander polynomial of a 3-manifold}\label{Alex_and_BNS}

Let $M$ be a compact, orientable 3-manifold whose boundary, if any, is
a union of tori.  Let $B_M = H_1(\tilde{M}, \Z)$ be the Alexander
invariant of $M$.  Regarding $B_M$ as a $\Z[\ab(\pi_1 M)]$ module, we
can form the BNS-invariant $\Sigma_{B_M}$ which I will denote by
$\Sigma_A$.  In Section~\ref{MoralReason}, I defined $\Sigma_A$ in
case where $H_1(M, \Z)$ is torsion free, and that definition was
slightly different.  In the torsion free case, $B_M = G'/G''$ where $G
= \pi_1(M)$.  Thus only difference between the two definitions is that
one is the BNS invariant with respect to $\ab(G)$ and the other $G$.
Since $B_M$ is abelian and $G'$ acts trivially on it, the two
definitions agree.

In this section I will prove Theorem~\ref{fullB_AandSigma_A} which
computes $\Sigma_A$ from the Alexander polynomial $\Delta_M$.  Before
stating Theorem~\ref{fullB_AandSigma_A}, I need to discuss the unit
ball $B_A$ in the Alexander norm.   

Consider the Newton polytope $\Newt(\Delta_M)$ in $H_1(M, \R)$.  The
Alexander norm on $H^1(M, \R)$ can be defined as
\[
\norm{\phi}_A = \sup\setdef{\phi(x - y)}{x, y \in \Newt(\Delta_M)}.
\]
A polytope $P$ is \emph{balanced} about $0$ if it is invariant under
$v \mapsto -v$.  More generally, $P$ is balanced about a point $p$ if
the translate of $P$ by $-p$ is balanced about $0$.  Since $M$ is a
3-manifold, $\Delta_M$ is symmetric \cite{Blanchfield},
\cite[4.5]{Turaev75}, and hence $\Newt(p)$ is balanced about some
point $z_0$.  Then
\[
 \norm\phi_A = \sup\setdef{2 \phi(x - z_0) }{ x \in \Newt(\Delta_M)}
\]
and the unit ball in $\norm{\ \cdot\ }_A$ is
\[
B_A = \setdef{\phi}{ \text{$\phi(x - z_0) \leq 1/2$ for all $x \in \Newt(\Delta_M)$}}.
\]
Fix a basis of $H_1(M, \R)$ and identify $H^1(M, \R)$ with
$H_1(M, \R)$ via the dual basis.  Then $B_A$ is, after scaling by a
factor of 2, the classical polytope dual of $\Newt(\Delta_M)$ about
$z_0$.

Duality of polytopes in an $n$-dimensional vector space exchanges
faces of dimension $i$ with faces of dimension $n - i - 1$ (for more on
polytope duals, see \cite{Broendsted}).  A vertex $v$ of $\Newt(\Delta_M)$
becomes the top-dimensional face 
\[
F_v = \setdef{\phi}{ \text{$\phi(x - z_0) \leq 1/2$ for all $x \in
    \Newt(\Delta_M)$ and $\phi(v - z_0) = 1/2$}}.
\]
I can now state the theorem that relates $\Sigma_A$ and $B_A$.
\begin{theorem}\label{fullB_AandSigma_A}
  Let $M$ be a compact, orientable 3-manifold whose boundary, if any,
  is a union of tori.  Let $F_1,\dots, F_n$ be the top-dimensional
  faces of $B_A$ whose corresponding vertices of $\Newt(\Delta_M)$
  have coefficient $\pm 1$.  Then $\Sigma_A$ is exactly the projection
  to $S(\ab(\pi_1M))$ of the interiors of the $F_i$.
\end{theorem}

\begin{proof}
Let $Q = \ab(\pi_1(M))$.  I will show:
\begin{lemma}\label{BNS_equality}
Let $M$ be as above.  Then $\Sigma_A = \Sigma_{\Z Q/ (\Delta_M) }$.
\end{lemma}
Let me now deduce the theorem assuming the lemma.  By
Theorem~\ref{computing_Sigma}, the components of $\Sigma_{\Z
  Q/(\Delta_M)} $ correspond to the vertices of $\Newt(\Delta_M)$
  whose coefficients are $\pm 1$.  Such a vertex $v$ corresponds to:
\[
C_v = \bigcap \setdef{H_{vw^{-1}}}{\text{$w$ is a vertex of $\Newt(p)$
      distinct from $v$}}, 
\]
where $H_q$ is the hemisphere $\setdef{\class{\chi} \in S(Q)}{\chi(q) > 0}$.
To prove the theorem it suffices to show $C_v$ is the same as the
projection into $S(Q)$ of the interior of the face $F_v$ of $B_A$
corresponding to $v$.  Translate $\Newt(\Delta_M)$ so it is balanced
about $0$---this doesn't change $C_v$ or $\norm{\ \cdot\ }_A$.  Now note
that the cone over the interior of $F_v$ is
\[
\setdef{\phi}{\text{$\phi(v) > \phi(w)$ for all vertices $w$ of
    $\Newt(p)$ distinct from $v$}}.
\]
It's easy to see that this cone projects to $C_v$ in $S(Q)$.  This
proves the theorem modulo the lemma.  Let's go back and prove the
lemma.

\begin{proof}[Proof of Lemma~\ref{BNS_equality}] The idea of the
proof is that Proposition~\ref{Alex_prop} says that $B_M$ is close, in
some sense, to $\Z Q/(\Delta_M)$.  Using the properties in
Section~\ref{BNS_definitions}, we have  (notation changed for clarity):
\[
\Sigma_A = \Sigma(B_M) = \Sigma\left(\Z Q/\Ann(B_M)\right) = 
 \Sigma\left(\Z Q\left/\sqrt{\Ann(B_M)}\right.\right).
\]
For any finitely generated module $B$ we have $\sqrt{\Ann(B)} =
\sqrt{E_0(B)}$, and so 
\[
\Sigma_A = \Sigma\left(\Z Q\left/\sqrt{E_0(B_M)}\right.\right).   
\]

Let $m$ be the augmentation ideal of $\Z Q$.  Since $\Z Q/m = \Z$, the
invariant $\Sigma_{\Z Q/m}$ is all of $S(Q)$.  So for any ideal $I$,
we have $\Sigma\left(\Z Q/(I \cap m)\right) = \Sigma\left(\Z Q/I\right) \cap
\Sigma\left(\Z Q / m\right) = \Sigma\left(\Z Q/I\right)$.  Thus
\[
\Sigma_A = \Sigma\left(\Z Q\left/\left(\sqrt{E_0(B_M)} \cap m\right)\right.\right).   
\]
By Proposition~\ref{Alex_prop}, $\sqrt{E_0(B_M)} \cap m = \sqrt{
  (\Delta_M ) } \cap m$, so
\[
\Sigma_A = \Sigma\left(\Z Q\left/ \left(\sqrt{(\Delta_M)} \cap m\right)\right. \right) =
\Sigma\left(\Z Q\left/ \left(\sqrt{(\Delta_M)}\right) \right. \right) =
 \Sigma\left(\Z Q/(\Delta_M)\right),
\]  as
  required.   This completes the proof of the lemma and thus the
  theorem.
\end{proof}
\renewcommand{\qed}{}
\end{proof}

\subsection{Comparison of $\Sigma_{G'}$ and
  $\Sigma_A$ when the homology is not free}\label{HomologyNotFree}

Let $M$ be a 3-manifold and $G$ its fundamental group.  In
Section~\ref{MoralReason}, I discussed the connection between
$\Sigma_{G'}$ and cohomology classes representing fibrations of $M$
over the circle.  This is true independent of whether $H_1(M, \Z)$ has
torsion.  In Section~\ref{MoralReason}, $\Sigma_{G'}$ and $\Sigma_A$
were compared under the assumption that $H_1(M, \Z)$ is free.  In this
case, it is easy to see $\Sigma_A \supset \Sigma_{G'}$, because $
\Sigma_A = \Sigma_{B_M}$, and $B_M = G'/G''$ is a quotient of $G'$.
When $H_1(M, \Z)$ is not free, the relation $\Sigma_A \supset
\Sigma_{G'}$ is still true, but not immediate since $B_M$ is a
quotient of the kernel of the map $G \to \ab(G)$, but that kernel
properly contains $G'$.

The purpose of this subsection is simply to prove the that $\Sigma_A
\supset \Sigma_{G'}$ for any $M$, and so show that the motivation given
in Section~\ref{MoralReason} makes sense regardless of whether $H_1(M,
\Z)$ is free.

\begin{proposition}
Let $M$ be a 3-manifold.  Then $\Sigma_A \supset \Sigma_{G'}$.
\end{proposition}

\begin{proof}
Let $N$ be the kernel of the map from $G$ to its \emph{free}
abelianization.  It is clear that $\Sigma_A \supset \Sigma_N$ as the
Alexander invariant $B_M$ is a quotient of $N$.  By Proposition~3.4 of
\cite{BieriNeumannStrebel87}, $\Sigma_N = \Sigma_{G'}$ and we are done.
\end{proof}

\section{Example of a specific link exterior in $S^3$}\label{SpecificExample}

\begin{figure}[thbp]
\centerline{\psfig{figure=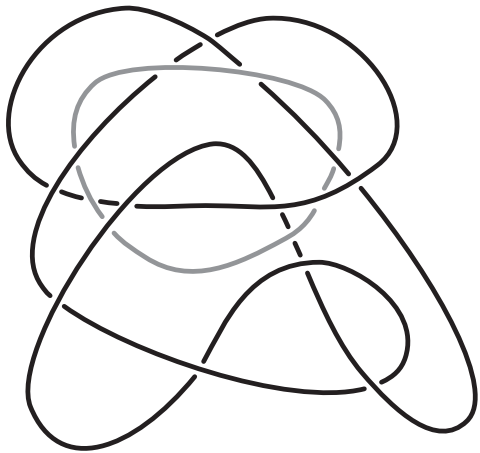, scale=\linkscale}}
\caption{\label{special_link}The link $L$ in $S^3$}
\end{figure}

Let $L$ be the link in Figure~\ref{special_link}.
Let $M = S^3 \setminus N(L)$ be the exterior of $L$.  In this section,
I'll show that $M$ is a fibered 3-manifold where the answer to \MyQ\ 
is no.  I will do this by explicitly computing the BNS invariants
$\Sigma$ and $\Sigma_A$, showing that $\Sigma$ is non-empty and that
each component of $\Sigma$ is properly contained in the corresponding
component of $\Sigma_A$.  The manifold $M$ is hyperbolic with volume
$8.997\ldots$, as can be checked with the program SnapPea
\cite{SnapPea}, or better, Snap \cite{Snap}, but I won't use this fact.
I found this example by a brute force search---the program SnapPea
was used to find many links whose fundamental groups have a
presentation with two generators and one relator.  For such groups, it is
easy to calculate $\Sigma$ and $\Sigma_A$ directly, as I will do
below, and eventually I came across this example.   

According to SnapPea, $\pi_1 M$ has a presentation with two generators
$a$ and $b$ and defining relation
\begin{equation*}
a^2 b a^{-1} b a^2 b a^{-1} b^{-3} a^{-1} b
  a^2 b a^{-1} b a b^{-1} a^{-2} b^{-1} a b^{-1} a^{-2} b^{-1} a b^3 a
  b^{-1} a^{-2} b^{-1} a b^{-1} a^{-1} b.
\end{equation*}
A meridian for the unknoted component is $b^{-1} a^{-1} b a^2 b
a^{-1} b a^2 b a^{-1} b^{-3}$ and a meridian for the other component
is $a^{-1} b^{-1}$.  

Let $X$ be the 2-dimensional CW-complex corresponding to the above
presentation.  Let $G = \pi_1(X)$.  The abelian group $\ab(G)$ is
freely generated by images of $a$ and $b$, and so $\Z[\ab(G)] =
\Z[\left<a, b\right>]$.

Let $\tilde{X}$ be the universal abelian cover of $X$.  It is natural
to think of the 1-skeleton of $\tilde{X}$ as the integer grid in
$H_1(X, \R)$.  Let $\delta$ be the lift of $a b a^{-1} b^{-1}$
starting at $0$, which freely generates the 1-chains of $\tilde{X}$ as
a $\Z[\ab(G)]$ module.  The 2-chains of $\tilde{X}$ are generated by
any lift of the 2-cell of $X$.  Let $\gamma$ be the lift of the
relator to 1-skeleton of $\tilde{X}$ starting at $0$, which is
homologous in the 1-skeleton to $(a^2 b - a b - a+1) \delta$.  Thus
\[
B_X = H_1(\tilde{X}, \Z) = \Z[\left<a,  b\right>] \big/ (a^2 b - a b - a + 1  ).
\]
So $\Delta_M = \Delta_X = a^2 b - a b - a+1$.  By
Theorem~\ref{fullB_AandSigma_A}, or, since $B_M$ cyclic,
Theorem~\ref{computing_Sigma} directly, we find that $\Sigma_A$ is all
of $S(\ab(G))$ except the four points $\left\{\pm [b^*], \pm [a^* -
  b^*]\right\}$, where $\{a^*, b^*\}$ is the dual basis to $\{a, b\}$.

\begin{figure}[bht]
  \centering 
  \parbox[b]{0.45\textwidth}{ 
    \centering
    \psfig{figure=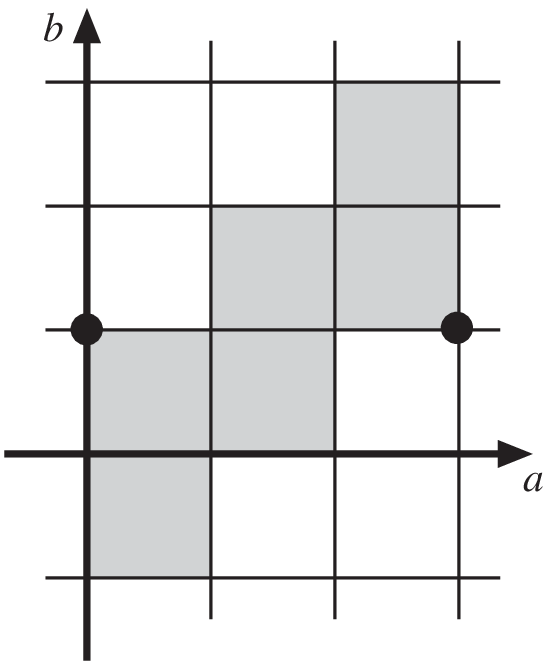, scale=\thirdpairscale}
   }\hspace{0.05\textwidth}%
   \parbox[b]{0.45\textwidth}{
     \centering
     \psfig{figure=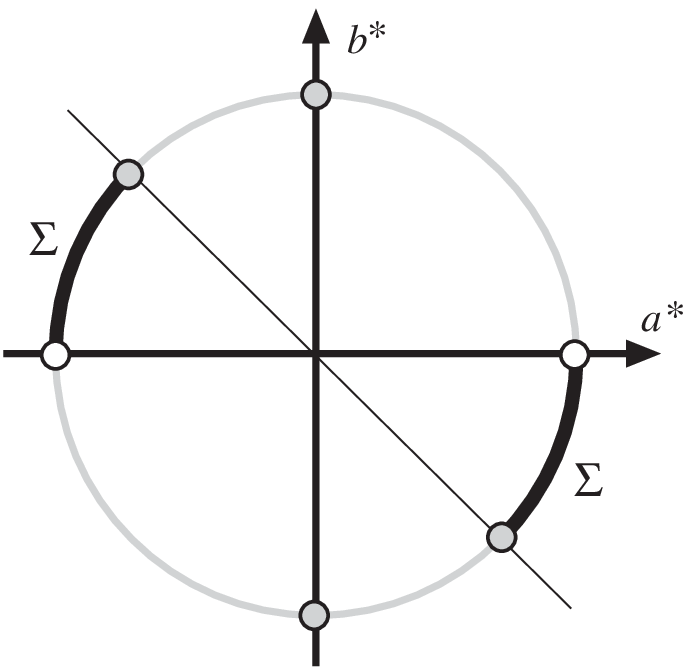, scale=\thirdpairscale}
   }

  \parbox[t]{0.45\textwidth}{ 
    \centering
    \caption{The region $C$.  The two dots are the 
      simple vertices $v_1$ and $v_2$.}\label{C}
   }\hspace{0.05\textwidth}%
   \parbox[t]{0.45\textwidth}{
     \centering
     \caption{$\Sigma \subset S(\ab(G))$
       consist of the two open intervals shown.  $\Sigma_A$ is the
       complement of the four grey dots.}\label{special_sigma}
   }
\end{figure}

To compute $\Sigma$, I'll use Brown's procedure for computing $\Sigma$
for any group with a 2-generator, 1-relator presentation \cite[\S
4]{Brown87}.  Think of the 1-skeleton of $\tilde{X}$ as the integer
grid in $H_1(X, \R)$.  Let $C$ be the convex hull of the $\gamma$, the
lift of the relator.  A vertex $v$ of $C$ is called simple if $\gamma$
passes through $v$ only once.  Figure~\ref{C} shows $C$ with the 2
simple vertices $v_1$ and $v_2$ marked.
Theorem~4.4 of \cite{Brown87} shows that in our case $\Sigma$
consists of two components $C_i$, for $i = 1, 2$, where
\[
C_i = \bigcap \setdef{H_{v_i w^{-1}}}{\text{$w$ is a vertex of $C$
      distinct from $v_i$}}.
\]
Thus $\Sigma$ is the union of the two open intervals pictured in
Figure~\ref{special_sigma}, and each component of $\Sigma$ is properly
contained in the corresponding component of $\Sigma_A$.  So $M$ shows
that the answer to \MyQ\ is no.

\small

\bibliographystyle{math} \bibliography{standard}

\end{document}